\documentclass[11pt]{article}

\usepackage{amsmath,amssymb}
\usepackage{amstext}
\usepackage{amsbsy}
\usepackage{graphics}
\usepackage{graphicx}
\usepackage{color}
\usepackage{epstopdf}
\usepackage{latexsym}
\usepackage{verbatim}

\newcommand{\sbar}{\bar{s}}
\newcommand{\etabar}{\bar{\eta}}
\newcommand{\fb}{\mathbf{f}}
\newcommand{\gb}{\mathbf{g}}
\newcommand{\pb}{\mathbf{p}}
\newcommand{\qb}{\mathbf{q}}
\newcommand{\bb}{\mathbf{b}}
\newcommand{\cb}{\mathbf{c}}
\newcommand{\db}{\mathbf{d}}
\newcommand{\gtilde}{\tilde{\mathbf{g}}}
\newcommand{\Rb}{\mathbf{R}}
\newcommand{\Rhat}{\mathbf{\hat{R}}}
\newcommand{\shat}{\mathbf{\hat{s}}}
\newcommand{\nhat}{\mathbf{\hat{n}}}
\newcommand{\xbar}{\bar{\mathbf{x}}}
\newcommand{\xb}{\mathbf{x}}
\newcommand{\K}{\operatorname{K}}
\renewcommand{\L}{\operatorname{L}}
\renewcommand{\S}{\operatorname{S}}

\newcommand{\ub}{\mathbf{u}}
\newcommand{\uinf}{\ub_{\infty}}
\newcommand{\real}{\mathbb{R}}
\newcommand{\eeps}{\varepsilon}

\title{Accurate evaluation of integrals in slender-body formulations
  for fibers in viscous flow}
\author{Anna-Karin Tornberg}
\date{}

\begin{document}
%\bibliographystyle{plain}
%\nocite{*}

%\pagestyle{myheadings}
%\thispagestyle{plain}
%\thispagestyle{empty}

\maketitle

\begin{abstract}
A non-local slender body approximation for slender flexible fibers in
Stokes flow can be derived, yielding an integral equation along the center lines
of the fibers that involves a slenderness parameter. 
The formulation contains a so-called finite part singular integral,
and can in the case of several fibers or evaluation of the flow field
require the evaluation of nearly singular integrals. 

We introduce a numerical technique to accurately and efficiently
evaluate the finite part integral. This technique can be applied
combined with any panel based quadrature rule and will add no
additional cost except for a small precomputation of modified quadrature
weights.  We also show how a related technique that was recently
introduced can be applied for the evaluation of the nearly singular
integrals.

\end{abstract}

%\begin{itemize}
%\item Use product integration for the special quadrature scheme? Or
%  semi-analytic?
%\end{itemize}

\section{Introduction}

Non-local slender body theory describes the
motion of flexible fibers or filaments in viscous flows. It is based
on Stokes equations, and can be used when the inertia of both fluid
and fibers can be neglected, i.e. for small Reynolds numbers. 
Slender body theory exploits the slenderness of the fibers, and is
more accurate the more slender the fiber is. 

Derivations can be found in \cite{Gotz00,KR76,Jo80}. The result is an
integral equation along the fiber center line with a slenderness
parameter $\eeps=r/L$, where $r$ is a representative radius of the
fiber and $L$ its length. Johnson \cite{Jo80} showed that this
equation is asymptotically accurate to $O(\eeps^2 \log \eeps)$
under some assumptions on the tapering
of the fiber towards the free ends.
Extending the equations to several fibers, the equations have an
asymptotic accuracy of $O(\eeps)$ when including both a Stokeslet and a
doublet kernel (the Laplacian of the Stokeslet) \cite{Gotz00}.

Shelley and Ueda \cite{SU96,SU00} were the first to construct a numerical method
based on a non-local slender body approximation. They did so for a
closed filament (i.e. with no free ends) with its motion constrained
to a plane in 3D space, studying the dynamics as the filament was set
to grow everywhere along its length. 
Tornberg and Shelley \cite{TS04}  extended this work to consider multiple
interacting slender fibers with free ends in a three dimensional
Stokes flow. A numerical method that included a semi-implicit treatment
in time was introduced, which eliminated the severe constraint on the
time step size that arise from the elasticity. 
Nazockdast et al. \cite{Nazock2017}, further improved on this discretization, in
both space and time. For a further and recent discussion on both
experiments, numerical methods and the applications of flexible
fibers in fluid see the recent review \cite{duRoure2019}.

In this paper, we will focus on the numerical evaluation of integrals
in the slender body formulation. In \cite{TS04}, Tornberg and Shelley
introduced a regularization of the so-called finite part integral in the
non-local operator to remove a solvability condition. This
regularization also makes the integral non-singular, but it is still
nearly singular, requiring care in its evaluation. In \cite{TS04}, a
piecewise linear approximation of the density was assumed between grid
points used to discretize the fiber such that integrals over each
subinterval could be evaluated analytically. The same approach was 
later used in \cite{Nazock2017}.
 
In this short note, we will focus on how to accurately evaluate the
original integral, leaving the option to regularize the equations
decoupled from the technique to accurately evaluate the
integral. Specifically for very small values of $\eeps$, the
regularization should not be needed.

We will introduce a method based on product integration to evaluate
the finite part integral. This technique was introduced by Helsing and Ojala
\cite{Helsing2008} to evaluate the harmonic single and double layer potentials in
2D. The technique was later extended to the Stokes equations in 2D by Ojala
and Tornberg \cite{Ojala2015}, and used in the simulation of viscous drops. 

Another integral that appears in the slender body integral formulation
contains the Stokeslet. It needs to be evaluated either to obtain the
fluid velocity in the field in a post-processing step, or already in
the solution step if multiple interacting fibers are considered. If the
fibers get close, this integral gets nearly singular and needs special
treatment for accurate evaluation. 
Here, we will use a method recently developed in \cite{afKB2019}, that is also
an extension of the technique introduced by Helsing and Ojala \cite{Helsing2008}. 

After giving some preliminaries, we will discuss the special
quadrature for the finite part integral in section \ref{sec:non-loc},
including its validation. Then we will turn to the nearly-singular
Stokeslet integral. For these near-singularities a special quadrature
method is developed in \cite{afKB2019}, and here we only describe and apply that
method to our case.

\section{Preliminaries}

Let the centerline of a fiber be parameterized by arclength $s \in
[0,L]$, where $L$ is the length of the fiber, and let
$\xb(s,t)=(x(s,t),y(s,t),z(s,t))$ describe the fiber centerline at time
$t$. 
We introduce the slenderness ratio $\varepsilon=r/L$, where $r$ is the
radius of the fiber. 
Given a background velocity $\uinf(\xb,t) \in \real^3$ of the fluid with viscosity
$\mu$, the non-local slender body
approximation \cite{Gotz00,KR76,Jo80}, gives the relation between the velocity
of the fiber centerline and the force per unit length $\fb(s,t) \in \real^3$, 
\begin{equation}
8 \pi \mu \left ( 
\frac{\partial \xb(\sbar,t)}{\partial t}
-\uinf(\xb(\sbar,t),t) \right)=-\Lambda[\fb](\sbar)-\K[\fb](\sbar),
\quad \sbar \in [0,L].
\label{eqn:SBeqn}
\end{equation}
The local operator $\Lambda[\fb](\sbar)$ is given by
\begin{equation}
\Lambda[\fb](\sbar)=
[-c(I+\shat \shat(\sbar))+
(2-\shat \shat(\sbar))] \fb(\sbar),
\label{eqn:SBLocOp}
\end{equation}
where $c=\log(\varepsilon^2 e)$ ($c<0$), $\shat(\sbar)$ is the unit tangent vector at $s=\sbar$, and
$\shat \shat$ is the dyadic product, i.e. $(\shat \shat)_{kl}=\shat_k \shat_l$.
The next term is a non-local operator, that is defined as 
\begin{equation}
K[\fb](\sbar)=\int_0^L \left[ 
\left(  \frac{I+\Rhat(s,\sbar) \Rhat(s,\sbar)}{|\Rb(s,\sbar)|} 
%%+ \frac{\Rb(s,\sbar) \Rb(s,\sbar)}{|\Rb(s,\sbar)|^3} 
 \right) \fb(s)
- \left( \frac{I+\shat \shat(\sbar)}{|s-\sbar|} 
\right) \fb(\sbar) 
\right] \,ds,
\label{eqn:SBNonLocOp}
\end{equation}
where we have introduced the notation 
$\Rb(s,\sbar)=\xb(s)-\xb(\sbar)$ and $\Rhat=\Rb/|\Rb|$. 
We are here suppressing the dependence on time in the notation. 

The non-local operator $\K$ introduced in (\ref{eqn:SBNonLocOp}) is a
so-called finite part integral. Each part is singular at $s=\sbar$ and
the integral is well defined only then the integrand is
kept as the difference between the two terms. 

The fluid velocity in a field point $\xbar$ can be approximated by 
\begin{equation}
8 \pi \mu \left ( 
\ub(\xbar)-\uinf(\xbar) \right)=-\int_0^L 
\left( 
\frac{I}{|\Rb(s)|} + \frac{\Rb \Rb(s)}{|\Rb(s)|^3}  \right) \fb (s) \, ds, 
\label{eqn:fluidvel}
\end{equation}
where $\Rb(s)=\xbar-\xb(s)$ and $\Rb \Rb(s)$ is again a dyadic
product. This kernel is the Stokeslet. Sometimes, a so-called Stokes
doublet is added with an $\varepsilon^2/2$ coefficient \cite{Gotz00,TS04}.

If we consider more than one fiber, this integral will give the
velocity contribution from one fiber at a point on another. In this
case, (\ref{eqn:SBeqn}) will be extended to a coupled system for all fibers
with this interaction term \cite{TS04}. 

This integral will be easy to resolve when $\xbar$ is far from the
fiber, but the integrand will get increasingly peaked as $\xbar$ moves
close to the fiber. This nearly singular case will appear e.g. when
two fibers are interacting at a close distance. 

\subsection{Regular quadrature}

For smooth integrals, we use a regular quadrature rule. More
specifically, we will use a composite Gauss-Legendre quadrature rule,
but other quadrature rules e.g. such as Clenshaw-Curtis quadrature based on 
Chebyshev polynomials that was used in \cite{Nazock2017} can be used. 

Consider the integral of a smooth function 
\[
I=\int_0^L \phi(s) \,ds.
\]
Let us now split the interval $[0, L]$ into $M$ intervals of equal size
$\Delta s=L/M$, and write the integral as
\begin{equation}
I=\sum_{m=1}^M \int_{(m-1)\Delta s}^{m \Delta s} \phi(s)  \,ds 
= \frac{\Delta s}{2}\sum_{m=1}^M \int_{-1}^{1} \phi(s(\eta))  
 \,d\eta.
\label{eqn:ivalsplit}
\end{equation}
On each panel, we have introduce a local parameter $\eta \in
[-1,1]$, such that $s_m(\eta)=(m-1)\Delta s + \frac{\Delta s}{2} (\eta +1)$ for $s
\in [m-1,m]\Delta s$, $m= 1,\ldots,M$. 
We now introduce the Gauss-Legendre quadrature nodes $\eta_{\ell}$
and weights $w_{\ell}$, $\ell=1,\ldots,n_{GL}$. In this note, we will
use a $16$ point Gauss-Legendre rule, i.e. $n_{GL}=16$, but other
orders can be used. 
With this, we approximate
\[
\int_{-1}^{1} \phi(s_m(\eta))  \,d\eta \approx \sum_{\ell=1}^{n_{GL}}
w_{\ell} \phi(s_m(\eta_{\ell})). 
\]

This quadrature rule will accurately approximate the integral as long
as the integrand is smooth and can be resolved with the underlying
discretization. 

\section{The non-local operator}
\label{sec:non-loc}

The finite part integral (\ref{eqn:SBNonLocOp}) can not be accurately
evaluated using a regular quadrature method. We will start by rewriting the
integral and then apply a method introduced by Helsing and Ojala
\cite{Helsing2008} for a semi-analytical treatment. 

\subsection{Rewriting the non-local operator}

Let us first consider a simpler operator defined by the integral 
\begin{equation}
\L[f](\sbar)=\int_0^L \frac{f(s)-f(\sbar)}{|s-\sbar|}
 \,ds,
\label{eqn:Lop0}
\end{equation}
Rewriting this integral as
\begin{equation}
\L[f](\sbar)=\int_0^L g_0(s,\sbar) \frac{s-\sbar}{|s-\sbar|}
 \,ds,
\label{eqn:Lop1}
\end{equation}
with 
\begin{equation}
g_0(s,\sbar)=\frac{f(s)-f(\sbar)}{s-\sbar}
\label{eqn:g0}
\end{equation}
it is easy to see that 
\begin{equation}
\lim_{s \rightarrow \sbar} g_0(s,\sbar)=f'(\sbar). 
\label{eqn:limit_g0}
\end{equation}

Similarly, we can rewrite the full $\K$ operator in (\ref{eqn:SBNonLocOp})
on this form, 
\begin{equation}
\K[\fb](\sbar)=\int_0^L \gb(s,\sbar) \frac{s-\sbar}{|s-\sbar|}
 \,ds,
\label{eqn:Kop1}
\end{equation}
where 
\begin{equation}
\gb(s,\sbar) = 
\left[
\left( I+\Rhat(s,\sbar) \Rhat(s,\sbar) \right) 
\frac{|s-\sbar|}{|\Rb(s,\sbar)|}
\fb(s)
-\left( I+\shat(\sbar) \shat(\sbar) \right) \fb(\sbar)
\right]
\frac{1}{s-\sbar}
\label{eqn:gdef}
\end{equation}

%Note that for a straight fiber, we can express $\gb(s,\sbar)=(I+\shat \shat(\sbar))\gb_0(s,\sbar)$
%with $\gb_0$ as $g_0$ in ($\ref{eqn:g0}$) generalized to a
%vector valued $\fb(s)$ in a straight forward manner. 

The limit $\lim_{s \rightarrow \sbar} \gb(s,\sbar)$ exists, and to
find it, we first subtract and add the term
$\left( I+\shat(\sbar) \shat(\sbar) \right) \fb(s)$ inside the square
bracket in ($\ref{eqn:gdef}$).  We then write
\[
\gb(s,\sbar) = \gb_1(s,\sbar) + \left( I+\shat(\sbar) \shat(\sbar)
\right) \gb_0(s,\sbar), 
\]
where $\gb_0$ is defined as  $g_0$ in ($\ref{eqn:g0}$) but with a
vector valued $\fb(s)$. Hence we have $\lim_{s \rightarrow \sbar}
\gb_0(s,\sbar)= \fb'(\sbar)$. 
There is a finite limit also for $\gb_1$, and to determine the limit, we
Taylor expand around $s=\sbar$. Adding the results together, we get
\begin{equation}
\lim_{s \rightarrow \sbar} \gb(s,\sbar)=
\frac{1}{2} ( \xb_s \xb_{ss}(\sbar) + \xb_{ss} \xb_{s}(\sbar)
)\fb(\sbar)
+(I+\shat \shat(\sbar)) \fb'(\sbar),
\label{eqn:limit_g}
\end{equation}
where subscripts $s$ denote derivatives with respect to arclength,
and hence $\xb_s=\shat$ and $\xb_{ss}=\kappa \nhat$, where $\kappa$ is
the curvature and $\nhat$ the principal normal. 

The operators in Eq (\ref{eqn:Lop1}) and Eq (\ref{eqn:Kop1})
are both now defined by an integral on the same form, with a smooth scalar or vector
valued function multiplying the kernel $(s-\sbar)/|s-\sbar|$.  Next,
we will consider how to accurately evaluate these integrals.

\subsection{Special quadrature method}

Let us now consider the evaluation of 
\begin{equation}
I_{\phi}(\sbar)=\int_0^L \phi(s,\sbar) \frac{s-\sbar}{|s-\sbar|}
 \,ds.
\label{eqn:Iphi}
\end{equation}
With $\phi=g_0$ as defined in (\ref{eqn:g0}), this defines the operator
$L$ in (\ref{eqn:Lop1}). If we instead let $\phi$ denote the $x$, $y$
or $z$ component of $\gb$ as defined in (\ref{eqn:gdef}), the integral
yields the corresponding component of the operator $K$ in 
(\ref{eqn:Kop1}). 

Dividing into subintervals as in (\ref{eqn:ivalsplit}), we have 
\[
I_{\phi}(\sbar)=\sum_{m=1}^M \int_{(m-1)\Delta s}^{m \Delta s} \phi(s,\sbar)  \frac{s-\sbar}{|s-\sbar|}
 \,ds =  
\frac{\Delta s}{2} \sum_{m=1}^M \int_{-1}^{1} \phi(s_m(\eta),s(\etabar))  \frac{\eta-\etabar}{|\eta-\etabar|}
 \,d\eta.
\]
Now, consider a panel with $s \in [m-1,m]\Delta s$. If the evaluation
point $\sbar$ lies outside of this interval,
$(s-\sbar)/|s-\sbar|$ will be constant over the full interval,
as there will be no shift in sign. 
Since $\phi(s,\sbar)$ is smooth, the full integrand will be smooth over
this interval, and regular quadrature can be used. 
For $\sbar \in [m-1,m]\Delta s$, the integrand has a discontinuity, and
we will use product integration for accurate results. Let 
\[
\etabar=-1+\frac{2}{\Delta s} \left(\sbar-(m-1)\Delta s\right),
\]
s.t. $\etabar \in [-1,1]$. 
For short, denote $\phi_{\etabar}(\eta)=\phi(s_m(\eta),s(\etabar))$ and consider the evaluation of 
\begin{equation}
I_m(\etabar)=\int_{-1}^{1} \phi_{\etabar}(\eta)  \frac{\eta-\etabar}{|\eta-\etabar|}
 \,d\eta.
\label{eq:Imdef}
\end{equation}
Expanding $\phi_{\etabar}(\eta)$ into a polynomial with $n_{GL}$ terms, 
\[
\phi_{\etabar}(\eta)=\sum_{k=0}^{n_{GL}-1} c_k \eta^k, \eta \in [-1, 1],
\]
we get 
\begin{equation} 
I_m(\etabar)= \sum_{k=0}^{n_{GL}-1} c_k \int_{-1}^{1} \eta^k  \frac{\eta-\etabar}{|\eta-\etabar|}
 \,d\eta = \sum_{k=0}^{n_{GL}-1} c_k q_k(\etabar),
\label{eq:ImdefCoeff}
\end{equation}
where
\begin{equation} 
q_k(\etabar) = \int_{-1}^{1} \eta^k
\frac{\eta-\etabar}{|\eta-\etabar|} \, d \eta =\frac{1+(-1)^{k+1} -2 \etabar^{k+1}}{k+1}.
%%q=(1+(-1).^(pvec+1)-2*tt(j).^(pvec+1))./(pvec+1);
\label{eq:qkdef}
\end{equation}

\subsection{Precomputation of modified quadrature weights}

We assume that we discretize the slender body integral equation using
a Nystr\"{o}m method. This means that we collocate the equation at the
quadrature nodes.  Hence, on panel $m$, we assume that we have $\phi$
evaluated at the Gauss-Legendre nodes $\eta_{\ell}$,
$\ell = 1, \ldots,n_{GL}$.

We define three column vectors, $i)$ $\pb$, $ii)$ $\cb$ and $iii)$
$\qb(\etabar)$, containing the values of $\{\phi(\eta_{\ell})\}_{i=1}^{n_{GL}}$, 
$\{c_k\}_{k=0}^{n_{GL-1}}$, and $\{q_k(\etabar)\}_{k=0}^{n_{GL-1}}$, respectively. 
The coefficients in $\cb$ are the solution to the Vandermonde system 
\[
A \cb = \pb, 
\]
where column number $k$ of the Vandermonde matrix $A$ contains the
values of $\{ (\eta_{\ell})^{k-1}\}_{\ell=1}^{n_{GL}}$.
With this, we continue from  Eq (\ref{eq:ImdefCoeff}) and write 
\[
I_m(\etabar)=\cb^T \qb(\etabar) = (A^{-1}\pb)^T \qb(\etabar)=\pb^T (A^{-T}\qb)=
\pb^T \bb(\etabar), 
\]
where in the last step, we have defined the vector $\bb$ as the
solution to 
\[
A^T \bb(\etabar) = \qb(\etabar). 
\]

When we solve the integral equation, $\etabar$ will in turn take the
values of all Gauss-Legendre quadrature nodes, since we collocate at
these nodes. Hence, for one reference panel, we can solve
$A^T \bb(\eta_{\ell}) = \qb(\eta_{\ell})$, $\ell=1,\ldots,n_{GL}$, 
to find the taget specific weights $\bb(\eta_{\ell})$ for each
Gauss-Legendre node. 
Then the integral $I_m(\eta_{\ell})$ in Eq (\ref{eq:Imdef}) can simply
be evaluated using these weights,
\begin{equation} 
I_m(\eta_{\ell})= \sum_{k=1}^{n_{GL}} b_k(\eta_{\ell}) \phi_{\eta_{\ell}}(\eta_k),
\label{eq:ImdefFinal}
\end{equation}
where $\phi_{\eta_{\ell}}(\eta_k)=\phi(s(\eta_k),s(\eta_{\ell}))$ as
introduced above (\ref{eq:Imdef}).

%On each panel, introduce a local parameter $\eta \in
%[-1,1]$, such that $s=(m-1)\Delta s + \frac{\Delta s}{2} (\eta +1)$ for $s
%\in [m-1,m]\Delta s$, $m= 1,\ldots,M$. 
%We now write
%\[
%I(\sbar)=\sum_{m=1}^M \int_{(m-1)\Delta s}^{m \Delta s} g(s,\sbar)  \frac{s-\sbar}{|s-\sbar|}
% \,ds 
%=
%\sum_{m=1}^M \int_{-1}^{1} g(s(\eta),s(\etabar))  \frac{\eta-\etabar}{|\eta-\etabar|}
% \,d\eta,
%\]

%ONLY NEED SPECIAL TREATMENT WHEN TAGRET PT INSIDE SUB-IVAL

Note that only once do we need to compute the $16$ target specific
weights $b_k(\eta_{\ell})$ for the $16$ target values
$\eta_{\ell}$. They can then be used to integrate over each panel for
any of the discrete target points within the panel. Remember, for
evaluation points outside of the panel, regular quadrature can be used.

\subsection{Validation and numerical tests}

G\"{o}tz \cite{Gotz00} has shown that the operator $\L$ in (\ref{eqn:Lop0})
diagonalizes under the Legendre polynomials $P_n$. Scaled to the interval
$s \in [0,L]$ this result yields
\begin{equation}
\L[\tilde{P}_n](\sbar)=-\lambda_n \tilde{P}_n(\sbar), \quad n=0,1,\ldots
\end{equation}
where $\tilde{P}_n(s)=P_n(-1+2s/L)$ and 
\[
\lambda_n=\lambda_{n-1}+\frac{2}{n}, \quad n>0 \quad \mbox{and} \quad \lambda_0=0.
\]

%G\"{o}tz [CITE] has shown that the operator 
%\begin{equation}
%\Lbar[\phi](\etabar)=\int_{-1}^1 \frac{\phi(\eta)-\phi(\etabar)}{|\eta-\etabar|}
% \,d \eta,
%\label{eqn:Lbar}
%\end{equation}
%diagonalizes under the Legendre polynomials $P_n$, i.e. 
%\[
%\Lbar[P_n](\etabar)=-\lambda_n P_n(\etabar). 
%\]
%For $\L$ in \ref{eqn:Lop0}, this yields, 
%\[
%\L[\tilde{P}_n](\sbar)=-\lambda_n \tilde{P}_n(\sbar),
%\]
%where $\tilde{P}_n(s)=P_n(-1+2s/L)$. 

We will start to investigate the performance of our special quadrature
on this example, since there is an exact result to compare to. 
Note that the shape of the fiber does not enter this integration. 
In all our tests, we use a $16$-point  Gauss-Legendre rule on each
panel, $n_{GL}=16$. 
\begin{figure}[htbp]
\centering
\resizebox{!}{4cm}{\includegraphics{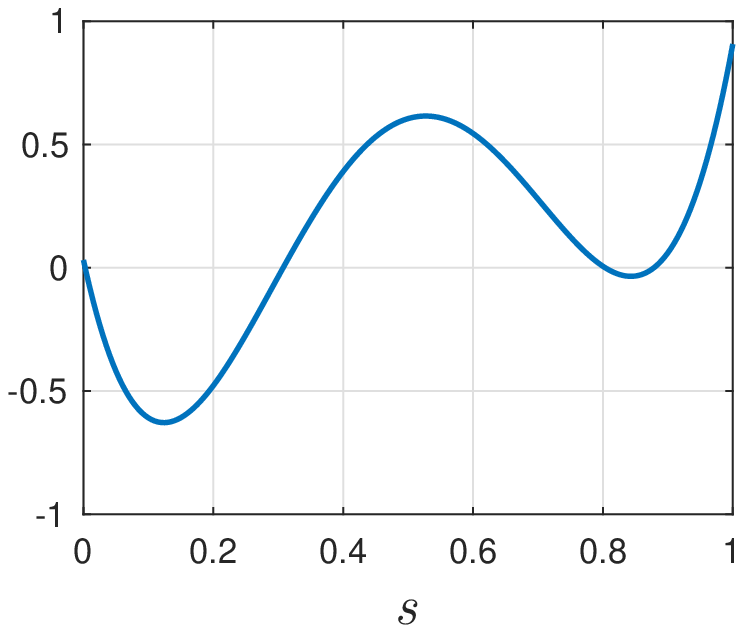}}
\resizebox{!}{4cm}{\includegraphics{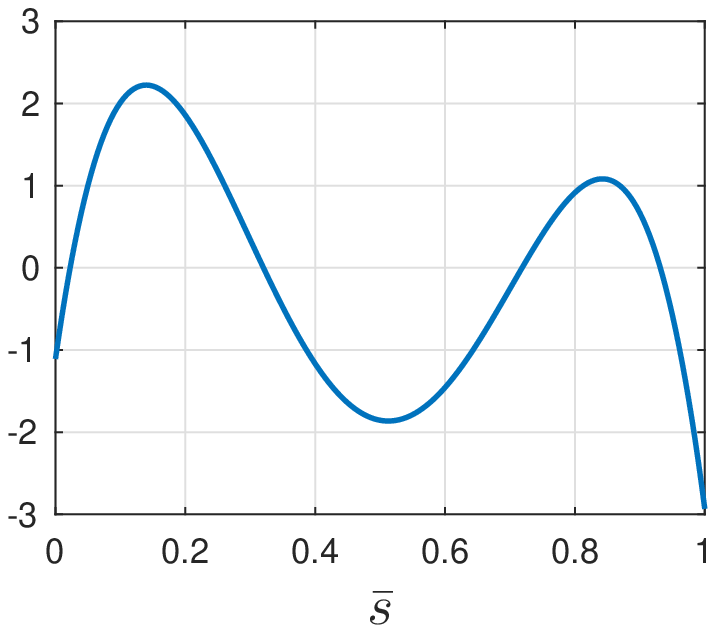}}
\caption{Left: $f(s)$ as in (\ref{eqn:Lexp_f}) with $P=5$. 
Right: $\L[f](\sbar)$ as in (\ref{eqn:Lop0}), plotted vs $\sbar$. 
}
\label{fig:f_Kbar_P5}
\end{figure}

In this first example, we set 
\begin{equation}
f(s)=\sum_{n=0}^{P-1} \alpha_n \tilde{P}_n(s), 
\label{eqn:Lexp_f}
\end{equation}
with the coefficients $\alpha_n$ random numbers between  $-1$ and $1$,
and we set $L=1$. 
See figure \ref{fig:f_Kbar_P5} for $f(s)$ and
$\L[f](\sbar)$. 

As expected, the errors in the results from the special quadrature
are at round off level for this case, see Table  \ref{tab:Lop}. It
does not matter here if we use $1$, $2$, $4$ or $8$ panels, since
already $1$ panel is sufficient to resolve this integral. 
%\begin{table}[htbp]
%\centering
%\begin{tabular}{|c|c|c|}
%\hline
% & Maximum error     & Maximum error \\
%No of panels & special quadrature & Matlabs built in \\
%\hline 
%%& & \\ 
%1& $1.78 \cdot 10 ^{-15}$  & $2.17 \cdot 10 ^{-10}$\\
%2&   $1.78 \cdot 10 ^{-15}$   &  $1.41 \cdot 10 ^{-8}$\\
%4 & $1.78 \cdot 10 ^{-15}$  & $3.60 \cdot 10 ^{-6}$ \\
%8& $2.22 \cdot 10 ^{-15}$   &  $9.97 \cdot 10 ^{-4}$ \\
%\hline 
%\end{tabular}
%\caption{Error in evaluation of $\L[f](\sbar)$ with $f$ as in Figure
%  \ref{fig:f_Kbar_P5}. Error measured by comparing to value obtained
%using diagonalization result. Maximum error taken over all
%$\sbar_{\ell}$, the $16 \times$ (No of panels) Gauss Legendre points. }
%\label{tab:Lop}
%\end{table}

\begin{table}[htbp]
\centering
\begin{tabular}{|c|c|}
\hline
 No of panels & Maximum error    \\
M& special quadrature \\
\hline 
%%& & \\ 
1& $1.78 \cdot 10 ^{-15}$  \\
2&   $1.78 \cdot 10 ^{-15}$   \\
4 & $1.78 \cdot 10 ^{-15}$  \\
8& $2.22 \cdot 10 ^{-15}$    \\
\hline 
\end{tabular}
\caption{Error in evaluation of $\L[f](\sbar)$ with $f$ as in Figure
  \ref{fig:f_Kbar_P5}. Error measured by comparing the values obtained
  using special quadrature to the diagonalization result. Maximum error taken over all
$\sbar_{\ell}$, the $16 M$ Gauss Legendre points. }
\label{tab:Lop}
\end{table}

As we turn to the full operator $\K$, we have no diagonalization
result or another analytical result that we can use for validation and
we opted to try the built in Matlab routine \texttt{integral} for
adaptive integration.  For the simpler operator $\L$, it works well if
we manually split the integration interval in two parts, $[0,\sbar]$
and $[\sbar,L]$. Without this manual split of the interval, errors
fluctuate by orders of magnitude for different values of $\sbar$.
For $\K$, such a split also improves the result, but the error levels
are higher. With this, we can validate our results down to an error level of 
about $10^{-9}$. 

%In using this routine, we ask for both a relative and absolute error of
%$10^{-10}$, in one integral over the full interval length.  We
%consider for the case discussed above, for which we can compare to the
%diagonalization result.  In Figure \ref{fig:MatErr_Kbar_P5}, we plot
%the error in the evaluation of $\L[f](\sbar)$ versus $\sbar$. Matlabs
%built-in command \texttt{integral} has been used with the same error
%tolerances for each value of $\sbar$, but the error varies widely by
%several orders of magnitude.
%\begin{figure}[htbp]
%\centering
%\resizebox{!}{4cm}{\includegraphics{Figures/MatlabErrUniKbarP5.eps}}
%\label{fig:MatErr_Kbar_P5}
%\caption{
%Error in evaluating $\L[f](\sbar)$ as in (\ref{fig:f_Kbar_P5}) for 
%$\sbar_{\ell}=\ell/N_{uni}$, $\ell=0,\ldots,N_{uni}$ using Matlabs
%built in adaptive integration.
%}
%\end{figure}
%%We conclude that we cannot use this routine for validation. 

Another option for validation is to compare results obtained with the
special quadrature rule for different number of panels. 
However, the discrete Gauss-Legendre points $\sbar_{\ell}$ do not
coincide for different number of panels. 
In order to compare the results, we interpolate all results
to a uniform grid, on each panel using the naturally defined Legendre
polynomials. 

Introduce a uniform grid with $N_u=400$ points. Compute a reference
solution $\K^{Ref}(\sbar)$ for $\K[\fb](\sbar)$ with $128$ panels, and interpolate the
result to the values of $\sbar$ of this uniform grid,
$\sbar_{\ell}=\ell L/N_u$, $\ell=0,\ldots,N_u$. 
Now compute an approximation of $\K[\fb](\sbar)$ with $M$ panels 
($\K^{M}(\sbar)$), interpolate to the uniform grid, and define 
\begin{equation}
e_M=\max_{0 \le \ell \le N_u} \| \K^M(\sbar_\ell) - \K^{Ref}(\sbar_\ell) \|_2.
\label{eq:errdef}
\end{equation}

\begin{figure}[htbp]
\centering
\resizebox{!}{4cm}{\includegraphics{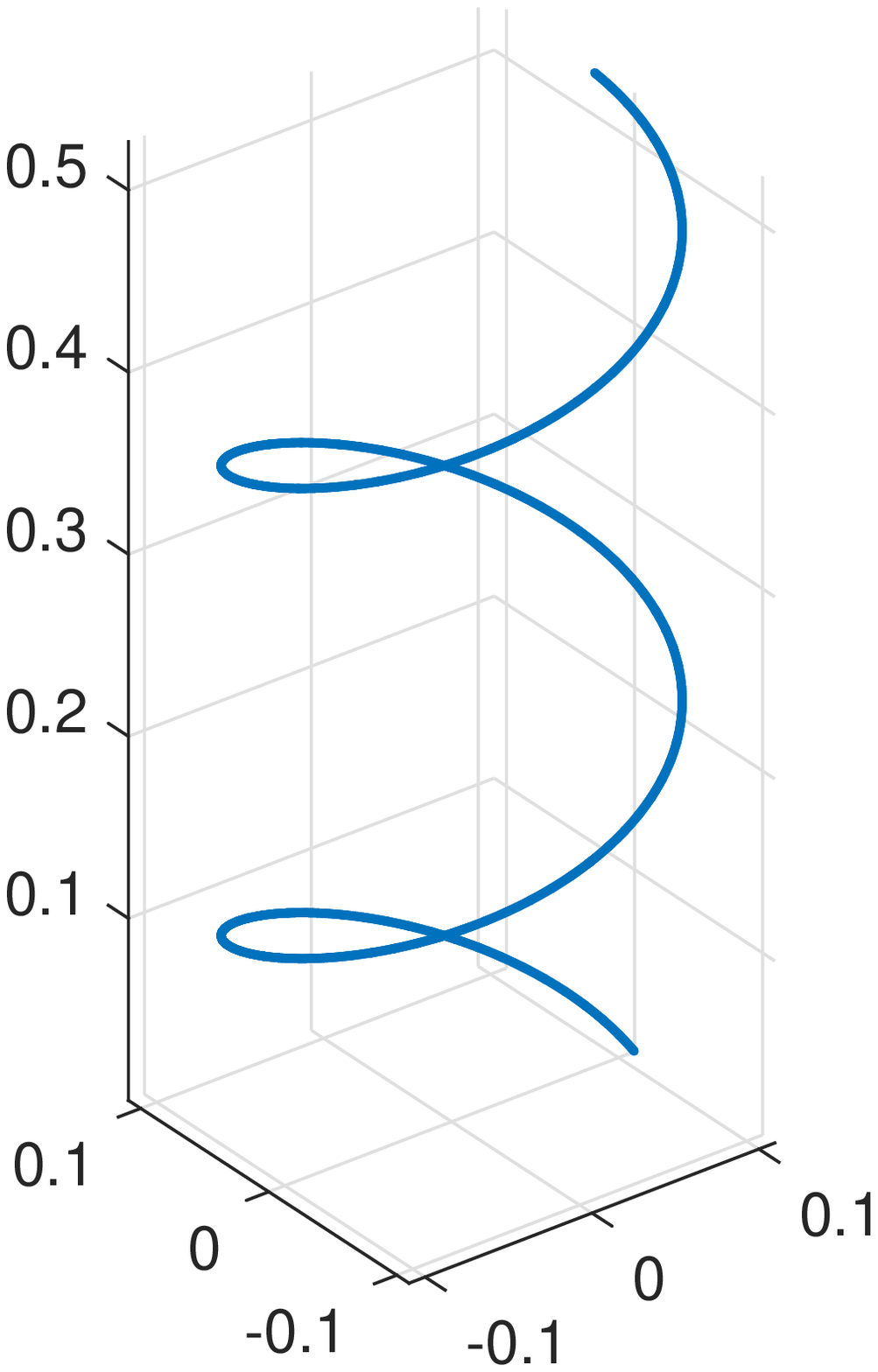}}
\resizebox{!}{4cm}{\includegraphics{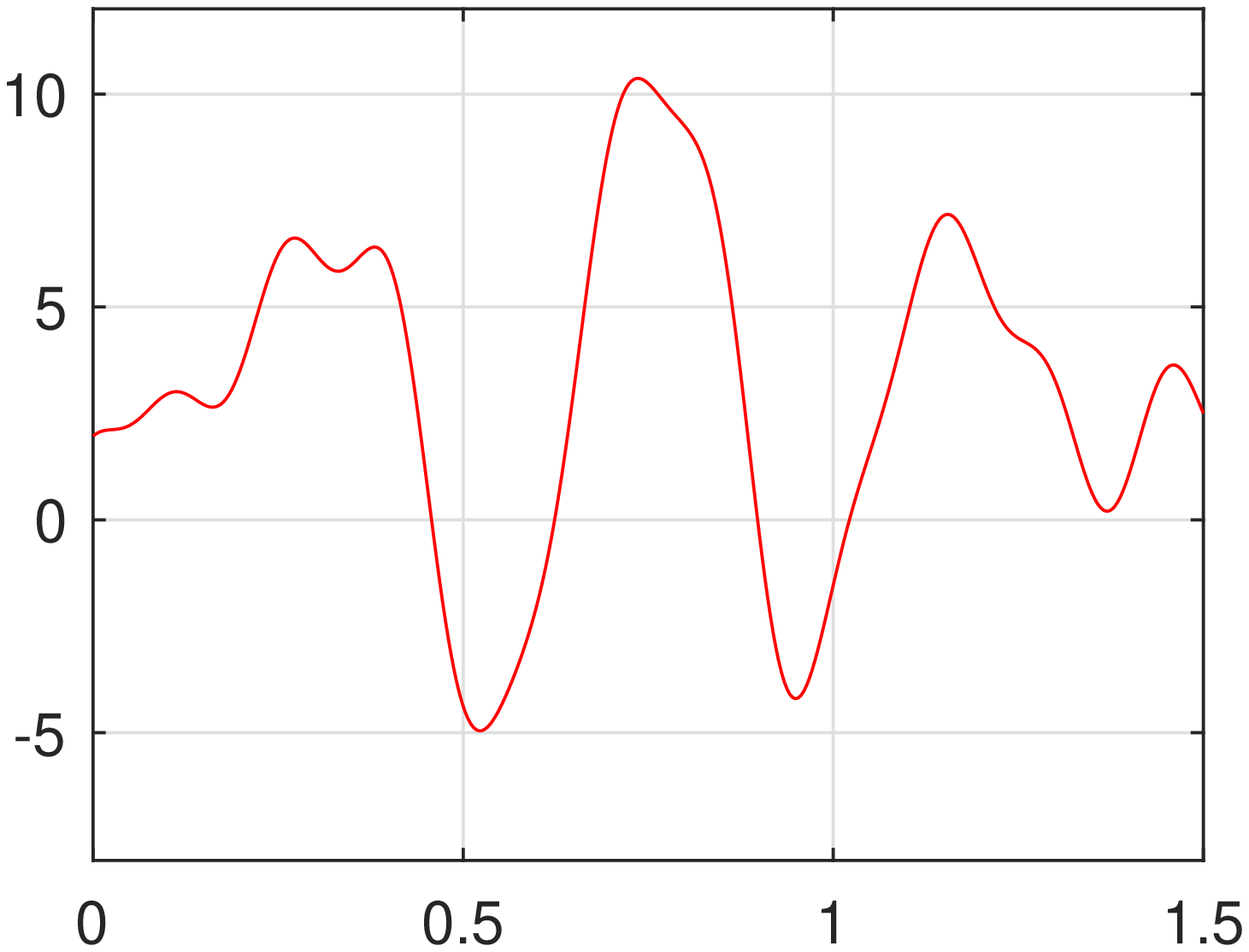}}
\resizebox{!}{4cm}{\includegraphics{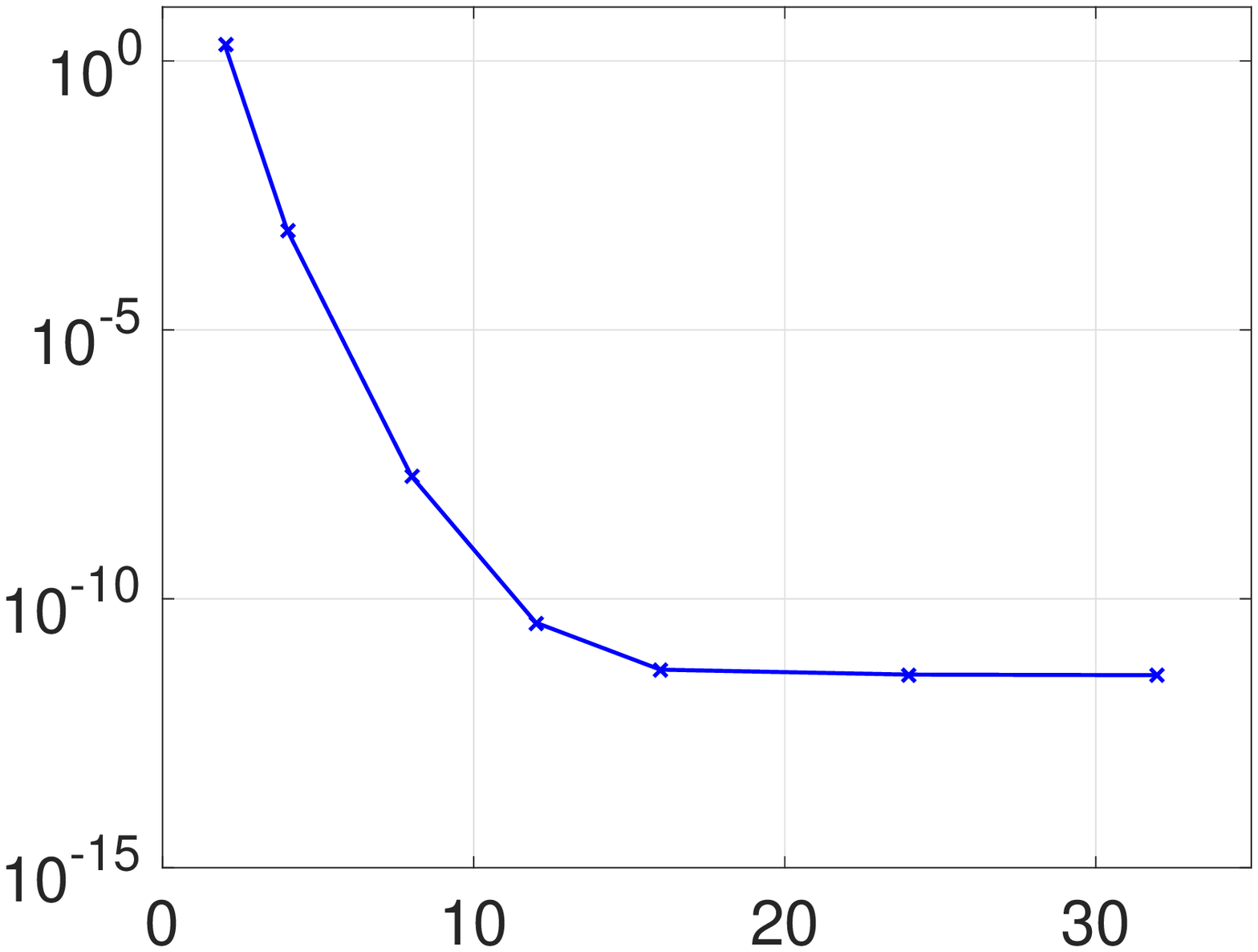}}
\caption{Left: A helix with constant curvature $\kappa=8$ and torsion
  $\tau=3$, total arclength $L=3/2$. 
Center: 
$\K[\fb](\sbar)$
plotted versus $\sbar$ for $\fb$ in (\ref{eq:testf}). 
Right: 
The error $e_M$ (\ref{eq:errdef}) in approximation of $\K[\fb](\sbar)$
versus number of panels ($M$). 
\label{fig:full_op_uni}
}
\end{figure}

In figure \ref{fig:full_op_uni}, the results from such a convergence
test is shown. Here, the fiber is set to be a helix with constant
curvature and torsion, and the force is given by
$\fb(s)=(f_1(s),f_2(s),f_3(s)$, where 
\begin{equation}
  f_1(s) = \cos(2\pi s)^2+e^{-s}+e^{-L+s}, \quad
  f_2(s) = \sin(4 \pi s)^2, \quad
  f_3(s)  =e^{-2s}.
\label{eq:testf}
\end{equation}
Here, we see the rapid decay of the error as we increase
the number of panels. We should however note that the error that we
are measuring is not only the quadrature error, but includes also the
error introduced from interpolation to the uniform grid. We believe
that it is the interpolation error that makes the error curve flatten
out at a level around $10^{-12}$. 

\section{Fluid velocity and interacting fibers}

The Stokeslet integral in (\ref{eqn:fluidvel}) gets very difficult to
resolve as the evaluation point $\xbar$ gets close to the fiber, and
it is not possible to deal with this difficulty simply by refinement. 
Doubling the number of panels, the width of the error region above a
certain tolerance will approximately be halved, but the maximum error
will not decrease. 

\subsection{Special quadrature method for nearly singular integral}
\label{sec:StokesSpecQuad}

In \cite{afKB2019}, the same underlying idea from \cite{Helsing2008} that we have
already used for the finite part integral has been extended to deal
with nearly singular line integrals.  

Let us write the integral in  (\ref{eqn:fluidvel}) as
\begin{equation}
\S[\fb](\xbar) = 
\int_0^L \left( \frac{\fb (s)}{|\xbar-\xb(s)|} + \frac{\left(
      (\xbar-\xb(s)) \cdot \fb(s) \right) (\xbar-\xb(s))}{|\xbar-\xb(s)|^3}
\right) \, ds.
\label{eq:StokesInt}
\end{equation}
Splitting the interval into panels in correspondance to
(\ref{eqn:ivalsplit}), we write for $p=1,3$
\[
\int_0^L \frac{\gb_p (s)}{|\xbar-\xb(s)|^p} \, ds=
\frac{\Delta s}{2}\sum_{m=1}^M \int_{-1}^{1} \frac{\gb_p(s_m(\eta))}{|\xbar-\xb(s_m(\eta))|^p}
 \,d\eta, 
\]
with $\gb_p(s)$, $p=1,3$ as understood from the two terms in (\ref{eq:StokesInt}).

Denote $\xb(\eta)=\xb(s_m(\eta))$, introduce $R^2(\eta)=|\xbar-\xb(\eta)|^2$
 and consider the integral 
\[
I_p(\xbar)=
%%\int_{-1}^{1} \frac{g(\eta)}{|\xbar-\xb(\eta)|^p} \,d\eta=
\int_{-1}^{1} \frac{g(\eta)}{\left(R^2(\eta)\right)^{p/2}} \,d\eta=
\int_{-1}^{1} \tilde{g}(\eta) \frac{1}{\left(\omega(\eta)\right)^{p/2}} \,d\eta
\]
where 
\[
%\tilde{g}(\eta) = g (\eta) \frac{\omega(\eta)}{\left(R^2(\eta)\right)^{p/2}}. 
\tilde{g}(\eta) = g (\eta) \left(\frac{\omega(\eta)}{R^2(\eta)}\right)^{p/2}. 
\]
The idea is similar to the finite part integral, to identify an
$\omega(\eta)$ such that $\tilde{g}(\eta)$ is regularized as compared
to the original integrand, expand $\tilde{g}(\eta)$ into a polynomial
and analytically evaluate the remaining integrals. 

Following \cite{afKB2019}, we define
\[
\omega(\eta)=(\eta-z_1)(\eta-\bar{z}_1)=|\eta-z_1|^2,
\]
where $\{z_1, \bar{z}_1\}$ is the complex conjugate root pair of
$R^2(\eta)$ that is closest to the interval $[-1,1]$.
The $q_k$:s corresponding to  (\ref{eq:qkdef}) will then be defined as
\[
q_k^p(z_1) = \int_{-1}^{1} \frac{\eta^k}{\omega(\eta)^{p/2}} \, d \eta =
\int_{-1}^{1} \frac{\eta^k}{|\eta-z_1|^p} \, d \eta, \quad p=1,3.
\]
Recursion formulas and a discussion about their numerical
evaluation is available in \cite{afKB2019}.

Hence, the structure is the same as before, but one needs to find
$z_1$. This is in \cite{afKB2019} done by using a Legendre expansion of each
component of $\xb(\eta)$ to define $R^2(\eta)$ combined with root
finding with Newton's method, see that paper for details. 

Note that the special quadrature is only needed when $\xbar$ is close
to the panel. Error estimates as derived in \cite{afKT2017,afKT2018} can be
used to estimate the error for the regular quadrature at any $\xbar$
and thereby determine when to switch to special quadrature. 

For the finite part integral, evaluation points are the Gauss-Legende
points on the panels. We can precompute the target specific weights
for one reference panel and use them for all panels. Here, we have no
prior knowledge of $\xbar$, and the special quadrature weights must be
computed as needed. This can however be done efficiently, see the
discussion in \cite{afKB2019}.

\subsection{Validation and numerical tests}

\begin{figure}[htbp]
\centering
\resizebox{!}{4cm}{\includegraphics{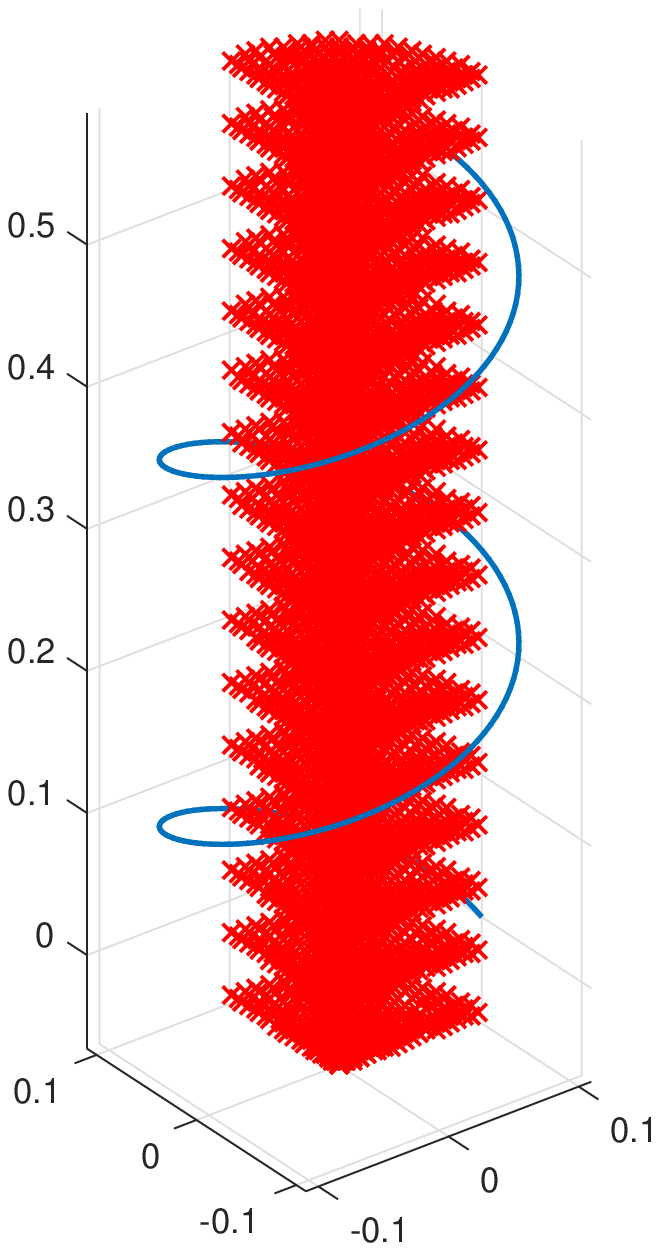}}
\resizebox{!}{4cm}{\includegraphics{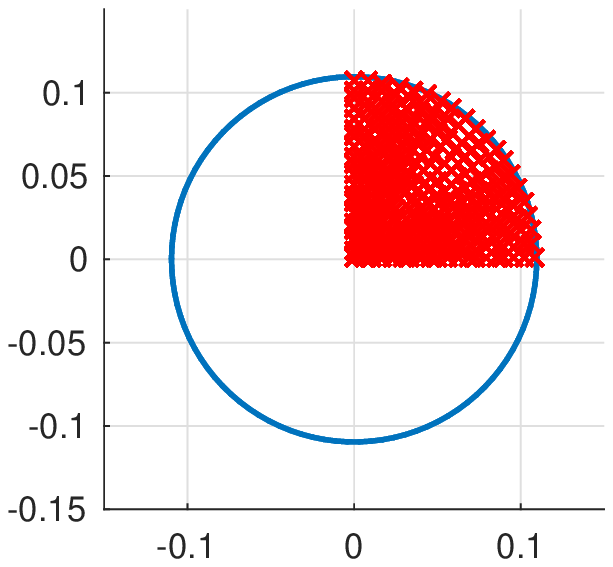}}
\caption{Left: The same helix as in figure \ref{fig:helix_field_points}, with field evaluation
  points. 
%%radius of the helix $8/73$.
$20 \times 20$ points are placed equidistant in polar coordinates
inside the projected circle of the helix (right) with closest distance
to boundary $2.2 \cdot 10^{-3}$. 
This is repeated for $16$ different $z$-values. 
Right: Projection onto the $xy$-plane.
}
\label{fig:helix_field_points}
\end{figure}

\begin{figure}[htbp]
\centering
\resizebox{!}{3.5cm}{\includegraphics{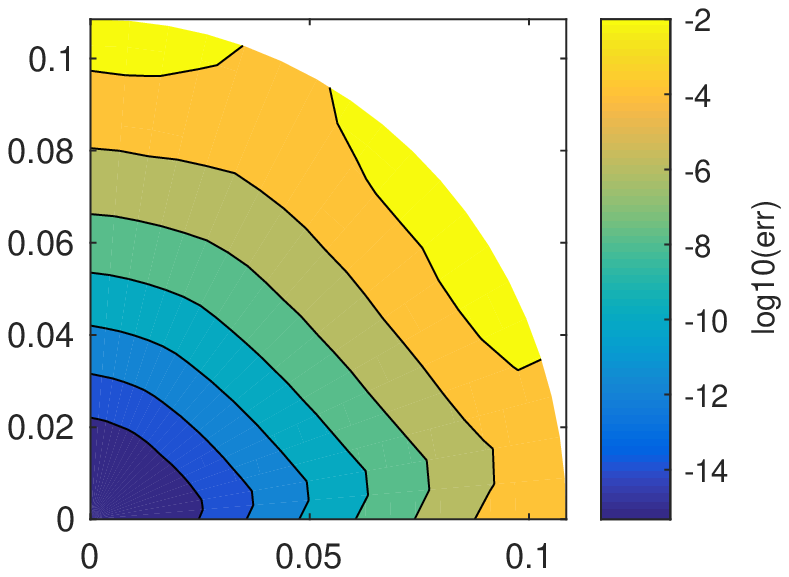}}
\resizebox{!}{3.5cm}{\includegraphics{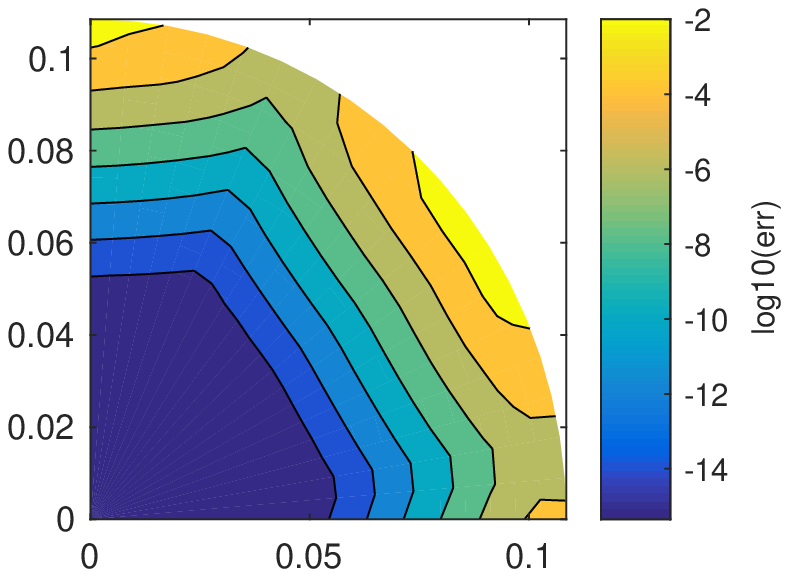}}
\resizebox{!}{3.5cm}{\includegraphics{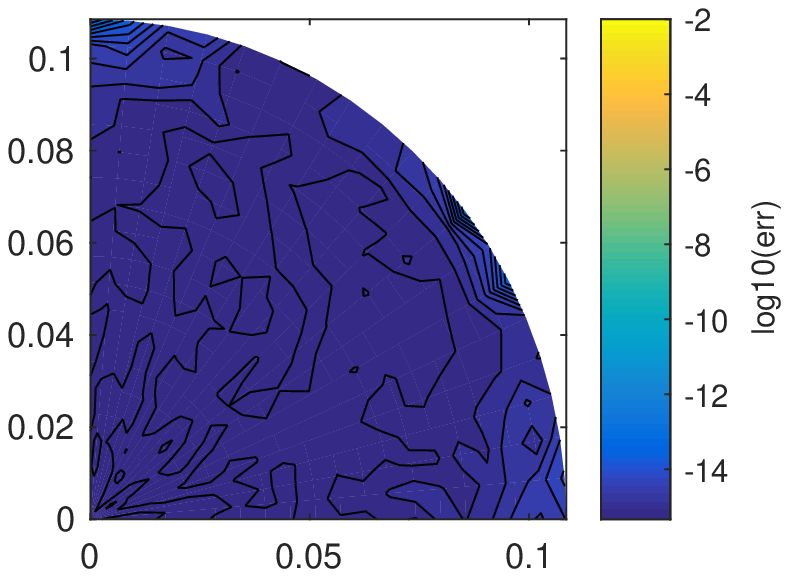}}
\caption{
Error in evaluation of $\S[\fb](\xbar)$: Regular Gauss-Legendre quadrature with $M=6$
panels (left), $M=12$ panels (middle) and special quadrature with $M=8$
panels (right). 
The error has been taken as the maximum over all $z$-values for each
evaluation coordinate $(x,y)$ as in figure
\ref{fig:helix_field_points} with $\fb(s)$ as in (\ref{eq:testf_simp}). }
\label{fig:field_errors}
\end{figure}

\begin{figure}[htbp]
\centering
\resizebox{!}{4cm}{\includegraphics{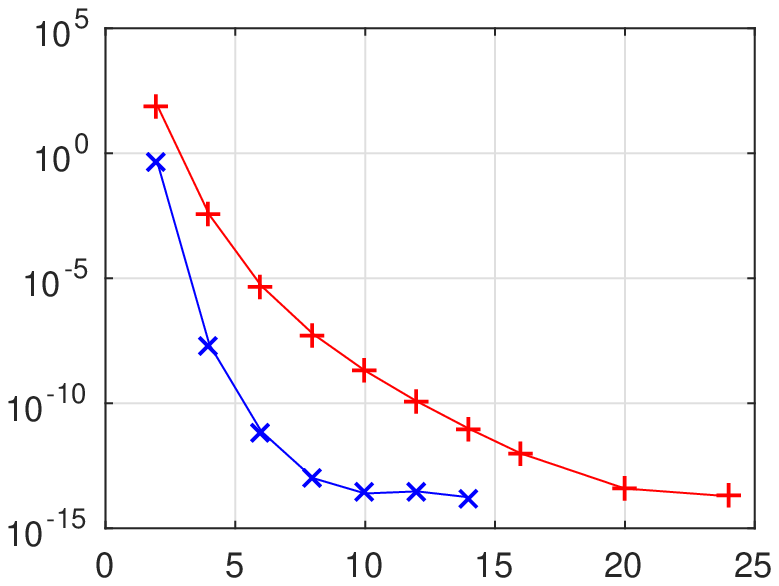}}
%\resizebox{!}{4cm}{\includegraphics{FiguresNearSing/err_vs_npan_field_kappa8.eps}}
%\resizebox{!}{4cm}{\includegraphics{FiguresNearSing/err_vs_npan_field_kappa8_compldensity.eps}}
\caption{ Maximum error in evaluation of $\S[\fb](\xbar)$ over all the
  evaluation points $\xbar$ in figure \ref{fig:helix_field_points}
  plotted versus the numer of panels $M$ for $\fb(s)$ as in
  (\ref{eq:testf_simp}) (blue line) and  (\ref{eq:testf}) (red
  line). 
% err_vs_npan_field_kappa8_2dens.eps, data in conv_data.mat
}
\label{fig:conv_off_surface}
\end{figure}

In the numerical examples, we use the same helix as in figure
\ref{fig:full_op_uni}. The helix is such that it projects onto a circle in the $xy$-plane. We
place evaluation points inside a quarter of that circle in different
$z$-planes, see figure \ref{fig:helix_field_points}.
We evaluate $S[\fb](\xbar)$ (\ref{eq:StokesInt}) for all these
evaluation points. We use $\fb(s)=(f_1(s),f_2(s),f_3(s)$, where 
\begin{equation}
  f_1(s) = x(s)+10 \quad
  f_2(s) = \sin(s), \quad
  f_3(s)  =\cos(s),
\label{eq:testf_simp}
\end{equation}
and also as defined in (\ref{eq:testf}), which is harder to
resolve. 
The integral (\ref{eq:StokesInt}) can be evaluated using Matlabs built
in adaptive quadrature (\texttt{integral}) to high precision, and the
error vector $\db=(d_1,d_2,d_3)$ is defined pointwise as the
difference to this reference solution. We define
$e(\xbar)=(d_1(\xbar)^2+d_2(\xbar)^2+d_3(\xbar)^2)^{1/2}$ and take
the maximum over all $z$-planes to display the errors in the
$xy$-plane in figure \ref{fig:field_errors}.
We take the maximum over all points to compute the maximum error
displayed in figure \ref{fig:conv_off_surface}.  From the results
using regular quadrature for $M=6$ and $M=12$ panels, we can see how
the error region shrinks, but how refinement fails to reduce the
errors closest to the boundary. The reason that the contours are not
circles is that the maximum is taken over a discrete set of
$z$-values. The special quadrature is keeping the errors very small
all the way up to the boundary.
In figure \ref{fig:conv_off_surface}, we can see the rapid decay of
the maximum error as we increase the number of panels with which we
discretize the helix. The special quadrature handles the singularity,
but the remaining density must also be well resolved for high
precision, and the figure shows the error for two different choices of
$\fb(s)$.

\section{Conclusions}
We have rewritten the integrand in the finite part integral as a
product of two factors. The first is now smooth and we have derived an
explicit formula for the limit at the problematic point. This factor is expanded as a polynomial, and the
integral over each term multiplied by the second factor can be
evaluated analytically. The coefficients in the polynomial expansion
are defined as the solution of a Vandermonde system. We show how we
can avoid to solve that system, and instead only once solve a sequence of transposed 
Vandermonde systems as a precomputation step, thereby defining
modified quadrature weights to be used in the numerical
evaluation. 
For a panel based quadrature with $n$ points on each panel, $n$ small
transposed Vandermode systems of size $n \times n$ must be solved, and
$n^2$ target specific quadrature weights must be stored.
Then the same modified quadrature weights can be used for all panels. 

The singularity is hereby treated analytically and does not cause any
error. The error from one panel will be determined by how well the first factor in
the integrand can be approximated by an $n-1$ degree polynomial over
that panel. In our numerical examples, we have used a $16$ point Gauss-Legendre
rule ($n=n_{GL}=16$). 

A technique that is similar in spirit but more complicated in its
detail is used for evaluation of the nearly singular integrals. Here,
the evaluation point can be any point close to the fiber, and no
precomputation is possible. Also in this case, we split the integrand into two factors,
where the first factor is to be approximated by a polynomial such that
the integral over each term can be analytically evaluated. 
To accomplish the split, one however needs to find a complex
conjugate pair by a root finding algorithm, and the analytical
evaluation makes use of recursion formulas that require care in their
numerical evaluation. This algorithm was introduced in \cite{afKB2019}. In
this paper, we have used it for the specific integrals that arise in
the slender body formulation, and have shown that it does indeed yield very
accurate results for evaluation points arbitrarily close to the
curve.

\section{Acknowledgements}

This work is dedicated to Professor Michael Shelley on the occasion of his 60th birthday.
The author wants to thank Ludvig af Klinteberg for sharing the 
implementation of the special quadrature method in \cite{afKB2019}.
This work was partially supported by the
G\"{o}ran Gustafsson Foundation for Research in Natural Sciences and Medicine, which is
gratefully acknowledged.

\bibliographystyle{plain}
\bibliography{slender}
\end{document}